\documentclass[12pt]{amsart}
\usepackage{amsfonts}
\usepackage{amsthm}
\usepackage{amssymb}
\usepackage{amsmath}
\usepackage[english]{babel}
\newtheorem{teo}{Theorem}
\newtheorem{prop}[teo]{Proposition}

\newtheorem*{teorI}{Theorem I}
\newtheorem*{teor0}{Theorem 1}

\newcommand{\fun}[5]{\renewcommand{\arraystretch}{1.5}
\begin{array}{crcl}
       #1: & #2 &\longrightarrow & #3 \\
      & #4 &\longmapsto & #5
       \end{array} \renewcommand{\arraystretch}{1} }

\newcommand{\co}{\mathbb{C}}

\newcommand{\re}{\mathbb{R}}

\newcommand{\qe}{\mathbb{Q}}

\newcommand{\cpt}[1]{\mathbb{C}P^{2}}

\newcommand{\pe}{\mathbb{P}}

\newcommand{\sing}{\mbox{Sing}}
\newcommand{\cl}[1]{\mathcal{#1}}
\newcommand{\ba}{\setminus}

\newcommand{\dr}{\mbox{$\partial$}}

\begin{document}
\title[Foliations in $\pe^n$ tangent to Levi-flat hypersurfaces]{On the dynamics of foliations in $\pe^n$ tangent to   Levi-flat hypersurfaces}
\author{Arturo Fern\' andez-P\' erez,  Rog\'erio Mol \& Rudy Rosas  }
\thanks{This work was supported by Pronex-FAPERJ, FAPEMIG, Universal CNPq and CAPES-Mathamsud. The first author named is partially supported by CNPq grant number 301635/2013-7.}
\date{\today}

\begin{abstract}
Let $\cl{F}$ be a codimension one holomorphic foliation in $\pe^{n}$, $n \geq 2$, leaving invariant a real analytic
 Levi-flat hypersurface $M$ with regular part $M^{*}$. Then every leaf of $\cl{F}$ outside $\overline{M^*}$ accumulates in $\overline{M^*}$.
\end{abstract}
\maketitle
\section{Introduction}
Let  $M \subset U \subset \co^{n}$ be a real analytic variety of real codimension one, where $U$ is an open set. Let $M^{*}$ denote its {\em regular} part, that is, the smooth part of $M$ of highest  dimension --- near each point $x \in M^{*}$, the variety $M$ is a manifold of real codimension one.
For each   $x \in M^{*}$, there is a unique complex hyperplane $\cl{L}_{x}$ contained in the tangent space $T_{x}M^{*}$.
This  defines a real analytic distribution $x \mapsto \cl{L}_{x}$ of complex hyperplanes in $TM^{*}$, known as the {\em Levi distribution}. When this distribution is integrable in the sense of Frobenius, we say that $M$ is a {\em Levi-flat} hypersurface. The resulting foliation, denoted by $\cl{L}$, is known as the {\em Levi foliation}. This concept goes back to E. Cartan, who proved that there are local holomorphic coordinates $(z_{1},\ldots,z_{n})$ around $x \in M^{*}$
such that $M^{*} = \{ {\rm Im}(z_{n}) = 0\}$ (\cite{cartan1933}, Th\'eor\`eme IV).  As a consequence, the leaves of the Levi foliation $\cl{L}$ have local equations $z_{n} = c$, for $c \in \re$. From the global viewpoint, they  are complex manifolds of codimension one immersed
in $U$.
Cartan's local trivialization provides an intrinsic way to extend the Levi foliation  to a non-singular holomorphic foliation  in a neighborhood of  $M^{*}$. Locally, we extend  $\cl{L}$ to a neighborhood of $x \in M^{*}$ as the foliation having, in the  coordinates $(z_{1},\ldots,z_{n})$, horizontal leaves $z_{n} = c$, for $c \in \co$.    Since $M^{*}$ has  real codimension 1, this   is the unique possible local extension of $\cl{L}$, so that these local extensions glue together yielding a foliation defined in whole neighborhood of $M^{*}$.
Nevertheless, it is  not true in general that $\cl{L}$ extends to a holomorphic foliation in a neighborhood of
$\overline{M^{*}}$, even if  singularities are admitted. There are examples of Levi-flat hypersurfaces whose Levi foliations extend to   $k$-webs in the ambient space (see \cite{brunella2007} and \cite{fernandez2012}).
When there is a holomorphic foliation $\cl{F}$ in $U$ which restricted to $M^{*}$ is the Levi foliation $\cl{L}$ we say either that $M$ is {\em invariant} by $\cl{F}$ or that $\cl{F}$ leaves $M$ {\em invariant}.

Local problems about singular Levi-flat hypersurfaces have been studied by many authors --- see for instance   \cite{burns1999}, \cite{fernandez2011}, \cite{fernandez2013},  \cite{lebl2013}  and the references within.
Germs of codimension one foliations at $(\co^{n},0)$ leaving invariant real-analytic   Levi-flat hypersurfaces  are well understood:  they are given by the levels of  meromorphic functions --- possibly holomorphic --- according to a  theorem by D. Cerveau and
A. Lins Neto   (see \cite{cerveau2011}
and also \cite{brunella2012}).
 In this note, our object is a globally defined singular holomorphic foliation $\cl{F}$ in $\pe^{n}$, $n \geq 2$, having an invariant real
analytic Levi-flat hypersurface $M$.  Our main result asserts that Levi-flat hypersurfaces in $\pe^{n}$ are attractors for the ambient foliation, in the sense that the leaves of $\cl{F}$ accumulates in $\overline{M^{*}}$. This is the content of:

\begin{teorI}
Let $\cl{F}$ be a complex foliation of codimension one in $\pe^{n}$, with $n \geq 2$, leaving invariant a real analytic Levi-flat hypersurface $M$.
Then every leaf of $\cl{F}$ accumulates in $\overline{M^{*}}$.
\end{teorI}

 We give a sketch of the proof of Theorem I, starting with the result in dimension two. First of all, by \cite{lebl2012}, the components of $\pe^{2} \ba \overline{M^*}$ are Stein Varieties, which, by a theorem of H\"ormander, are properly embedded in the affine space $\co^{5}$.
If there existed a leaf $L$ of $\cl{F}$ whose  closure does not
intersect $\overline{M^*}$, then $\overline{L}$ should contain a singular point of $\cl{F}$, otherwise $\overline{L}$ would yield a minimal set contained in a Stein variety, which   is not allowed.
A leaf that accumulates
in a singular point either contains a separatrix in its closure or is contained in a nodal separator (see \cite{camacho2013}). This allows us to use  the standard  Maximum Modulus Principle or its version for nodal separators
contained in Proposition \ref{nodalmaximum} to reach a contradiction with the compactness of $\overline{L}$. Finally, the $n$-dimensional result is obtained from the two-dimensional version by
restricting the foliation to a generic two-dimensional plane.

\section{Preliminaries}

 A holomorphic foliation of codimension one and degree $d$ in $\pe^{n}$ is induced, in homogeneous coordinates $(Z_{0}:Z_{1}:\cdots:Z_{n})$,
  by a 1-form $\omega$ whose coefficients are homogeneous polynomials of degree $d+1$ satisfying the following conditions:
\smallskip \par (i)  $\omega \wedge d \omega = 0$  (integrability);
\smallskip \par (ii)  $i_{\bf{r}} \omega \ =\ \sum_{i=0}^{n} X_{i} A_{i}(X) = 0$,
where ${\bf r} = X_{0} \dr / \dr X_{0} + \cdots X_{n} \dr / \dr X_{n}$ is the radial vector field (Euler's condition);
\smallskip \par (iii)  ${\rm codim} \, \sing(\cl{F}) \geq 2$,
 \smallskip \par \noindent where $\sing(\cl{F}) = \{A_{0} = A_{1} = \cdots = A_{n} = 0\}$ is the {\em singular set} of  $\cl{F}$. This means that $\omega$ defines, outside $\sing(\cl{F})$, a regular foliation of codimension one in $\co^{n+1}$ whose leaves are tangent
 to the distribution of tangent spaces given by $\omega$. Euler's condition assures that this foliation goes down to $\pe^{n}$.

  Let $M$ be an irreducible singular real-analytic Levi-flat hypersurface  in the complex projective space $\pe^{n}$.
Let $M ^{*}$ be  the   regular part of $M$, that is, the set of points near which $M$ is a nonsingular real-analytic hypersurface of real codimension one. We denote by   $\sing(M)$ the singular points  of $M$, points near which $M$ is not a real-analytic submanifold of any dimension. Note that, in general, $M^{*} \cup \sing(M) \subsetneq M$.

In our  approach to global Levi-flat hypersurfaces in $\pe^{n}$, an important tool is the use of geometric and analytic properties of Stein manifolds. Actually, we have the following result (see \cite{lebl2012}), which will play an essential
role the proof of Theorem I:
\begin{teor0}
\label{stein}
Let $M\subset \pe^{n}$ be a real-analytic Levi-flat hypersurface. Suppose that for every $p\in\overline{M^*}$ there exists a neighborhood $U_{p}$ and a meromorphic function $F_{p}$ defined in $U_{p}$ which is constant along the leaves of $M^*$. Then all the connected components of $\pe^n \setminus \overline{M^*}$ are Stein.
\end{teor0}
This result follows from a theorem of Takeuchi which asserts that
  an open set $U \subset \pe^{n}$ which is pseudoconvex is Stein
(see \cite{takeuchi1967}). The hypothesis on the existence of local meromorphic first integrals in the ambient assures that, at every point $p \in \overline{M^*}$, there exists a germ of complex
hypervariety contained in $\overline{M^*}$. This implies  pseudoconvexity for all connected components of $\pe^n \setminus \overline{M^*}$. Notice that, following Cerveau-Lins Neto's Theorem, such a condition is naturally
fulfilled when $M$ is tangent to a global foliation in $\pe^{n}$. Conversely,
the existence of local meromorphic first integrals allows a natural extension of the Levi foliation  to a neighborhood of $\overline{M^{*}}$ and a subsequent extension to the whole $\pe^{n}$, since  $\pe^{n} \ba \overline{M^{*}}$ has components which are Stein (see \cite{linsneto1999}).

A strong motivation to the theory of Levi-flat hypersurfaces is the study of minimal sets for foliations
(see \cite{camacho1988}, \cite{cerveau1993}).
If $\cl{F}$ is a singular foliation of dimension $r$ in a complex manifold $X$ of dimension $n>r$, then a compact non-empty subset
$\cl{M} \subset X$ is said to be a {\em minimal set} for $\cl{F}$ if the following properties are satisfied:
 \smallskip \par (i) $\cl{M}$ is invariant by $\cl{F}$;
 \smallskip \par (ii) $\cl{M} \cap \sing(\cl{F}) = \emptyset$;
 \smallskip \par (iii) $\cl{M}$ is minimal with respect to these properties.

 \smallskip \par Notice that if $L$ is a leaf of $\cl{F}$
such that $\overline{L}$ is compact and  $\overline{L} \cap \sing(\cl{F}) = \emptyset$, then $\overline{L}$
contains a minimal set for $\cl{F}$.
If $X$ is a Stein manifold, then $\cl{F}$ contains no minimal sets.
Indeed, a Stein manifold admits a $C^{\infty}$ strictly plurisubharmonic function $\phi$. If a minimal set
$\cl{M}$ existed, then the restriction of  $\phi$ to $\cl{M}$ would assume a maximum value at a point
$p \in \cl{M}$. If $L_{p}$ is the leaf of $\cl{F}$ containing $p$, the maximum principle for the plurisubharmonic
function $\phi_{|L_{p}}$ would force $\phi$  to be constant over $L_{p}$,  contradicting its
strict plurisubharmonicity. This fact has the following consequence:

\begin{prop}
\label{singintersection}
Let $\cl{F}$ be a complex foliation   in $\pe^{n}$  leaving invariant a real analytic Levi-flat hypersurface $M$.
Let $L$ be a leaf of $\cl{F}$. Then either $\overline{L} \cap \overline{M^{*}} \neq \emptyset$ or
$\overline{L} \cap \sing(\cl{F}) \neq \emptyset$.
\end{prop}
\begin{proof}
If neither of the alternatives were true, then $\overline{L}$ would be a minimal set contained in a Stein variety.
This is not allowed, as we commented above.
\end{proof}

For codimension one foliations, the existence of a non-singular Levi-flat hypersurfaces would imply the existence of a minimal set.
In $\pe^{n}$, for $n \geq 3$, there are neither non-singular real-analytic Levi-flat hypersurfaces, nor minimal sets, as proved in \cite{linsneto1999}. In dimension two, however, the existence of both real-analytic Levi-flats and minimal sets are
so far  open problems.

\section{Nodal separators}

A singular point $p$ for a local foliation $\cl{F}$ at at
$(\co^2,p)$ is {\em simple} if, given a vector field ${\bf v}$ that
induces $\cl{F}$ around $p$, the linear part $D{\bf v}(p)$ of ${\bf
v}$ at $p$ has eigenvalues $\lambda_1$ and $\lambda_2$ satisfying
one of the two possibilities:
\par  \smallskip (i) $\lambda_1 \neq 0$ and $\lambda_2 = 0$;
\par  \smallskip (ii) $\lambda_1 \lambda_2 \neq 0$ and $\lambda_1 / \lambda_2 \not \in \qe^+$.
\par \smallskip \noindent
Model number  (i) above is called {\em saddle node}, whereas a singularity satisfying (ii) is said to be {\em non-degenerate}.  Seidenberg's
Desingularization Theorem asserts that there is a finite sequence of
punctual blow-ups $\pi: (M,E) \to (\co^2,p)$, where $E = \pi^{-1}(0)$ is a normal
crossings divisor of projective lines and $M$ is a germ of complex surface around $E$, for which $\pi^{*}{\cl{F}}$, the strict transform of $\cl{F}$,
 is a foliation in $M$ whose singularities on $E$ are all  simple.  We say that
 a germ of complex foliation $\cl{F}$ with isolated singularity at $(\co^2,0)$ is a {\em generalized curve}
if there are no saddle-nodes in its desingularization (see \cite{camacho1984}).

If $p \in Sing(\cl{F})$ is a saddle-node, then there exists a smooth invariant curve $S$ invariant by $\cl{F}$
containing $p$, corresponding to the non-zero eigenvalue, the so-called strong separatrix --- to the zero eigenvalue
is associated a formal, possibly non-convergent, invariant curve named {\em weak separatrix}. In a neighborhood of $S \ba \{p\}$,
all leaves accumulate in $S$, so that they are not closed in $U \ba \{p\}$ for some neighborhood $U$ of $p$. This is incompatible, for instance, with the
existence of a local holomorphic first integral for $\cl{F}$, since, in this case, all leaves near $p$ would be closed
(see, for instance \cite{mattei1980}). Evidently, a germ of foliation $\cl{F}$ at $(\co^2,p)$ which admits a meromorphic first integral is a generalized curve, since its desingularization produces simple singularities admitting holomorphic first integrals. It follows from Ceveau-Lins Neto's Theorem that  a germ of foliation $\cl{F}$ at $(\co^2,p)$ leaving invariant
a germ of Levi-flat hypersurface is  a generalized curve.

Now, suppose that  $\cl{F}$ is a local foliation with a simple singularity at $p \in \co^{2}$
having non-zero eigenvalues $\lambda_{1}, \lambda_{2}$ such that $\lambda = \lambda_{2}/\lambda_{1} \in \re^{+} \ba \qe^{+}$.
Such a singularity is called {\em node}.
The condition on the eigenvalues puts this  singularity  in  the Poincaré domain, so that   $\cl{F}$ may be linearized by
holomorphic coordinates $(x,y)$, being  defined by the linear 1-form
\[\omega =  - \lambda y dx + x dy  .\]
The multivalued function $f_{\lambda} = x^{-\lambda} y$ is a first integral for $\cl{F}$. Outside the separatrices $x =0$ and $y = 0$, the function $|f_{\lambda}|$ is real analytic, so that $S_{c}: |f_{\lambda}| = c$, for a fixed  $c \in \re_{+}$, is  hypersurface, real analytic outside the separatrices, which is foliated by the leaves of $\cl{F}$, all of them dense in $S_{c}$.
 This set has the property that its complement in a small neighborhood of $p$ minus the local separatrices is not connected. It thus works as a barrier, preventing  local leaves in one component to pass to the other.
The set $S_{c}$ is  called {\em nodal separator}, following the terminology of D. Mar\'\i n and J.-F. Mattei (see \cite{marin2012}). The image of a nodal separator by a desingularization map is also called
a nodal separator, so that this concept  extends to non-simple singularities.

\par Sets with separation properties such as  nodal separators  do not exist for  non-nodal simple non-degenerate singularities.
Actually, when $\cl{F}$ is simple non-degenerate   at $p \in \co^{2}$,
then $\cl{F}$ has two separatrices. If, besides, $p \in \co^{2}$ is non-nodal, then the union of one of the separatrices with the {\em saturation} of a small  complex disc $\Sigma$ transversal to the other one --- that is, the union of all leaves of $\cl{F}$ intersecting $\Sigma$ --- is a neighborhood of $p$. Intuitively,   leaves that are sufficiently near one of the separatrices approach the other one.
Invariant objects with separation properties such as   nodal separators also fail to exist when
$p \in \co^{2}$ is a saddle node.

The main result in \cite{camacho2013} asserts that
if $\cl{F}$ is a germ of holomorphic foliation at $(\co^{2},p)$ with  isolated singularity
 and $\cl{I} \neq \{p\}$ is a closed invariant subset such that $p \in \cl{I}$, then $\cl{I}$ contains either a separatrix
or a nodal separator at $p$. In other words,  if $L$ is  a local leaf of $\cl{F}$ such that $p \in \overline{L}$, take
 $\cl{I} =  \overline{L}$ in order to conclude that $\overline{L}$ contains either a   separatrix or a nodal separator  at $p \in \co^{2}$.

Both separatrices and nodal separators are object of maximum modulus properties which will be useful
to our purposes. In the case of a separatrix, we have the following:
if a germ of analytic function $f \in \cl{O}_{0}$ and a germ of analytic curve
$S$ at $(\co^{2},0)$ are such that $|f|_{S}$ has a local maximum at some point $p \in S$ then $f$ is constant over $S$.
This follows from the standard Maximum Modulus Principle, after
possibly desingularizing  $S$.
Next, we propose  a version of the Maximum Modulus Principle for the universe of nodal separators.

\begin{prop}
\label{nodalmaximum}
 Let $S$ be  a nodal separator for a local foliation $\cl{F}$ in $(\co^{2},0)$. Let
$f \in \cl{O}_{0}$ be a holomorphic germ of function such that $|f|_{|S}$ has a maximum at $0 \in \co^{2}$. Then $f$ is constant.
\end{prop}
\begin{proof} By taking a desingularization, we may assume that $0 \in \co^{2}$ is a simple singularity, with
local coordinates $(x,y)$ such that $S$ has the  equation $|y|  =  |x|^{\lambda} $.
We only need to prove that that $|f|$ is constant over $S$. Actually, if this is so, then $f$ will be constant on each leaf
of $\cl{F}$ contained in $S$.  Each of these leaves accumulates in $0 \in \co^{2}$, so that $f \equiv f(0)$ on $S$. This shows that
the analytic set  $f = f(0)$ contains the non-analytic set $S$, which implies that $f \equiv f(0)$.

 Let $U$ be a neighborhood of $0 \in \co^{2}$, with $S \cap U$ connected, such that $|f(0)| \geq |f(p)|$ for all $ p \in S \cap U$. Notice that if there existed $p \in S \cap U$ with $p \neq 0$ such that  $|f(p)| = |f(0)|$, then $f$ would be constant on the leaf $L_{p}$  of $\cl{F}$ containing $p$. This leaf accumulates in $0$, so that $f \equiv f(0)$ on $L_{p}$. Again, the analytic set $f = f(0)$  would contain the non-analytic set $L_{p} \cup \{0\}$, which would imply $f \equiv f(0)$. Thus, we can suppose  that
$|f(0)| > |f(p)|$ for all $ p \in S \cap U$. By the
continuity of $f$, it is possible to obtain    a neighborhood $V$ of $S \ba \{0\}$, $V \subset U$, such that $|f(0)| > |f(p)|$ for all $p \in V$.

Let us fix $\epsilon > 0$ and a closed annulus $A = \{x \in \co; \rho_{1} \leq |x| \leq \rho_{2}\}$, for some $0 < \rho_{1} < \rho_{2}$,   such that $(x,y)$ in $V$ whenever
$x \in A$ and $| |y| - |x|^{\lambda}| < \epsilon$.
Thus, if $p/q \in \qe_{+}$ is sufficiently near
$\lambda$, then all the determinations of $(x,y^{p/q})$
will lie in $V$ for all $x \in A$. Let $S_{p/q}$ be the analytic curve of equation
$y^{q} = x^{q}$. We have that $0 \in S_{p/q}$ and, over the annulus $A$, $S_{p/q} \subset V$. Therefore, the  maximum of $f$ over $S_{p/q}$, for $|x| \leq \rho_{2}$, is reached at  some point $(x_{0},y_{0}) \in S_{p/q}$ with $|x_{0}| < \rho_{1}$. Now, the Maximum Modulus Principle applied to the analytic curve $S_{p/q}$ gives that $f \equiv f(0)$  over $S_{p/q}$. Finally,  curves such as $S_{p/q}$ are dense in the nodal separator, which allows us to
conclude that  $f$ is constant over $S$.
\end{proof}

\section{Proof of Theorem I}

We are now ready to proof  Theorem I:

\begin{teorI}
Let $\cl{F}$ be a complex foliation of codimension one in $\pe^{n}$, with $n \geq 2$, leaving invariant a real analytic Levi-flat hypersurface $M$.
Then every leaf of $\cl{F}$ accumulates in $\overline{M^{*}}$.
\end{teorI}
\begin{proof}
We first suppose that $\cl{F}$ is a foliation in $\pe^{2}$.
Let us suppose, by contradiction, that there exists a leaf $L$ such that $\overline{L} \cap \overline{M^{*}} = \emptyset$.
Then $\overline{L}$ is a compact set contained in a connected component $W$ of $\pe^{2} \ba \overline{M^*}$, which is a two-dimensional Stein variety. Notice that, by Proposition
\ref{singintersection}, necessarily $\overline{L}$ intersects $\sing(\cl{F})$.
A theorem of H\"ormander assures the existence of a proper holomorphic embedding $\phi:W \to \co^{5}$ (see \cite{hormander1990}). The restriction  of   $\cl{F}$ to $W$ is carried by $\phi$ to a foliation in the  closed surface $\phi(W) \subset \co^{5}$.
For simplicity, we use the same notation for $\cl{F}$ and its leaves
in $W \subset \pe^{2}$ and for their images in $\phi(W) \subset \co^{5}$.
Take a non-zero ${\bf v} \in \co^{5}$ and define
\[ \fun{f = f_{\bf v}}{\co^{5}}{\co}{z}{<z,{\bf v}>,} \]
where $<z,{\bf v}> = z_{1} v_{1} + \cdots + z_{5} v_{5}$,  $z = (z_{1},\ldots,z_{5})$ and
${\bf v} = (v_{1},\ldots,v_{5})$.
Since $\overline{L} \subset \co^{5}$ is compact, there exists $p \in \overline{L}$ where
$|f|_{\overline{L}}$ reaches its maximum value.

\smallskip
\noindent
\underline{\bf Assertion}: There exists a leaf $L_{1}$ of $\cl{F}$ with  $L_{1} \subset \overline{L}$, such that
$f_{|L_{1}}$ is constant.

\smallskip
\noindent We consider two cases:

\smallskip
\noindent
\underline{\bf 1st case}: $p$ is a regular point for $\cl{F}$. Take $L_{1}$ the leaf of $\cl{F}$ containing $p$. Since $L_{1} \subset \overline{L}$, the function $|f|_{|L_{1}}$ has a maximum value at $p$. Thus, $f_{|L_{1}}$ is constant by the Maximum
Value Principle.

\smallskip
\noindent
\underline{\bf 2nd case}: $p \in \sing(\cl{F})$. In this case,  by \cite{camacho2013},
$\overline{L}$ contains
 either a   separatrix or a nodal separator for $\cl{F}$ at $p$. Take $L_{1}$ to be the separatrix, in the first case, or
one of the leaves contained in the nodal separator, in the second. In both cases, $f_{|L_{1}}$ will be  constant
by the Maximum Value Principle (Proposition \ref{nodalmaximum}, for the nodal separator)

\smallskip The leaf $L_{1}$ found above is such that $f$ is constant over $\overline{L_{1}}$.
Thus, $\overline{L_{1}} \subset f^{-1}(f(p))$, which is a complex hyperplane in $\co^{5}$, that is,
an affine space isomorphic to $\co^{4}$.
Repeating  this  procedure three more times, we will  eventually find a leaf $L_{0} \subset \overline{L}$ such that $\overline{L_{0}} \subset \co$. This  gives a contradiction with the compactness of $\overline{L_{0}} $.

Now, suppose that $\cl{F}$ is a foliation in $\pe^{n}$, with $n \geq 3$, leaving invariant a real-analytic
Levi-flat hypersurface $M$. Let $L$ be a leaf outside $\overline{M^{*}}$. By \cite{camacho1992}, we can choose
a linear embedding $i: \pe^{2} \hookrightarrow \pe^{n}$ transversal to $\cl{F}$ and to $L$
 --- that is, $i^{*} \cl{F}$ is a foliation with isolated singularities having $i^{-1}(L)$ as a leaf.
Notice that $i^{-1}(M)$ is a real-analytic Levi-flat hypersurface whose Levi foliation is $i^{*} \cl{L}$, where $\cl{L}$ is the Levi foliation on $M$.  Now, we apply  the two-dimensional case in order to conclude that   $i^{-1}(L)$ intersects
  the closure of the regular part of $i^{-1}(M)$. This give straight
that $\overline{L} \cap \overline{M^{*}} \neq \emptyset$.
\end{proof}

\vskip 0.2 in


{\footnotesize
\medskip \medskip \medskip
\noindent
Arturo Fern\'andez-P\'erez  \\
Departamento de Matem\'atica \\
Universidade Federal de Minas Gerais \\
Av. Ant\^onio Carlos, 6627  \  C.P. 702  \\
30123-970  --
Belo Horizonte -- MG,
BRAZIL \\
arturofp@mat.ufmg.br

\medskip \medskip
\noindent
Rog\'erio  Mol  \\
Departamento de Matem\'atica \\
Universidade Federal de Minas Gerais \\
Av. Ant\^onio Carlos, 6627  \  C.P. 702  \\
30123-970  --
Belo Horizonte -- MG,
BRAZIL \\
rsmol@mat.ufmg.br

\medskip \medskip
\noindent
Rudy Rosas \\
Pontificia Universidad Cat\'olica del Per\'u\\
Av. Universitaria 1801, Lima, Per\'u\\
Instituto de Matem\'atica y Ciencias Afines IMCA\\
Jr. Los Bi\'ologos 245, Lima, Per\'u\\
rudy.rosas@pucp.edu.pe
}

\end{document}